\documentclass[conference]{IEEEtran}
\IEEEoverridecommandlockouts
\usepackage{cite,url}
\usepackage{amsmath,amssymb,amsfonts,mathtools}
\usepackage{algorithmic}
\usepackage{graphicx}
\usepackage{textcomp}
\usepackage{xcolor}
\usepackage{subfig, float}
\def\BibTeX{{\rm B\kern-.05em{\sc i\kern-.025em b}\kern-.08em
    T\kern-.1667em\lower.7ex\hbox{E}\kern-.125emX}}
\begin{document}

\title{Flattening the Duck Curve: A Case for Distributed Decision Making\\

\thanks{This material is based upon work partially supported by the MathWorks Mechanical Engineering Fellowship, U.S. Department of Energy under Award Number DE-IA0000025, NSF CPS-DFG Joint award 1932406, and the DIBRIS-UNIGE on institutional funding SEED 2019.}
}

\author{\IEEEauthorblockN{Rabab Haider\IEEEauthorrefmark{1},
Giulio Ferro\IEEEauthorrefmark{2},
Michela Robba\IEEEauthorrefmark{2}, and
Anuradha M. Annaswamy\IEEEauthorrefmark{1}}
\IEEEauthorblockA{\IEEEauthorrefmark{1}Mechanical Engineering, Massachusetts Institute of Technology, Cambridge, USA}
 \IEEEauthorrefmark{2} DIBRIS, University of Genoa, Genoa, Italy\\
Email: rhaider@mit.edu, giulio.ferro@edu.unige.it, michela.robba@unige.it, aanna@mit.edu
}


\maketitle

\begin{abstract}
The large penetration of renewable resources has resulted in rapidly changing net loads, resulting in the characteristic ``duck curve''. The resulting ramping requirements of bulk system resources is an operational challenge. To address this, we propose a distributed optimization framework within which distributed resources located in the distribution grid are coordinated to provide support to the bulk system. We model the power flow of the multi-phase unbalanced distribution grid using a Current Injection (CI) approach, which leverages McCormick Envelope based convex relaxation to render a linear model. We then solve this CI-OPF with an accelerated Proximal Atomic Coordination (PAC) which employs Nesterov type acceleration, termed NST-PAC. We evaluate our distributed approach against a local approach, on a case study of San Francisco, California, using a modified IEEE-34 node network and under a high penetration of solar PV, flexible loads, and battery units. Our distributed approach reduced the ramping requirements of bulk system generators by up to 23\%.
\end{abstract}


\section{Introduction}
The push towards decarbonization of the electric grid has seen an explosive growth in renewable energy installation. 
The state of California generates approximately 25\% of its electricity demand from solar resources, resulting in the characteristic ``duck curve'' in net system load \cite{duck}. This new operating condition, where dispatchable bulk resources must quickly meet the large and rapid change in electricity demand, introduces challenges to grid operators. Storage offers one solution, as a flexible and low-inertia resource, however remains too expensive for utility-scale system-wide adoption. Instead, we propose system operators look towards the distribution grid to provide some support. Small-scale consumer owned distributed energy resources (DERs) which include rooftop PV with inverters, demand response, and batteries can provide support to the bulk grid. 

Recent works in this area are \cite{sheha2020solving} which proposes a game theoretic framework for dynamic pricing, and \cite{calero2021duck} which shifts thermal cooling loads of residential units to mitigate duck curve effects. Where these works rely on local approaches with assets locally providing the mitigation mechanisms, our proposed work takes a systems-approach. We leverage distributed optimization to deploy resources over a larger spacial area of the distribution grid to enhance mitigation strategies while suitably accommodating global network constraints. 


The main contribution of the paper pertains to the distributed coordination of DERs across a large region of the distribution grid to mitigate the duck curve. This coordination is accomplished by leveraging the variable power factor setting of PV inverters, flexible loads to reduce consumption, and distributed storage devices such as community batteries. This central contribution is realized by formulating the global optimization problem as an optimal power flow, and employing a current-injection (CI) based linear model describing the power physics of the unbalanced distribution grid. The CI-based OPF is solved in a distributed fashion, using a Proximal Atomic Coordination (PAC) approach \cite{RomvaryTAC} with Nesterov acceleration, called NST-PAC. We evaluate the performance on a modified IEEE-34 node network, using load and generation data for San Francisco. Our results show the successful mitigation of the duck curve through the coordination of spatially distributed DERs.

In Sections~\ref{sec:CI} and \ref{sec:PAC} we introduce the CI power flow model for the unbalanced distribution grid, and the accelerated distributed optimization algorithm called NST-PAC. In Section~\ref{sec:case} we carry out a numerical case study, comparing local and distributed paradigms. We provide conclusions in Section~\ref{sec:conclusion}.


\section{Current Injection Model} \label{sec:CI}
The distribution grid is a highly unbalanced network, with many single- and two-phase lines and corresponding loads. As DER penetration increases, many of which are added to single-phase lines, the network will become more unbalanced. In light of this, power flow models which assume balanced flow \cite{gan2014exact} or those which are valid for only a small range of angle imbalances \cite{robbins2015optimal} cannot suitably describe the future grid. Further, the faster timescales of operation due to renewable variability and flexibility requirements need a model which is simple and results in OPF problems which are computationally tractable. For this reason, there is interest in developing scalable linear models for the unbalanced grid, and leveraging distributed computation and algorithms to solve the OPF problem. 

In this work, we utilize the Current Injection (CI) model, a linear model for multi-phase unbalanced distribution grids. The CI model takes a similar approach to the Bus Injection Model, wherein all nodal power injections are modeled as current injections, and all phasor variables are represented in Cartesian coordinates. The 3-phase impedance matrix is used to describe the self and mutual inductance between phases to model the coupling of phases that are common to a distribution grid. The non-convexity of the AC-OPF and the subsequent nonlinearity of SOCP and SDP convexification strategies typically used are avoided in the CI approach by leveraging McCormick Envelope (MCE) based convex relaxation \cite{mccormick_1976} for the bilinear power relations. The MCE uses the convex hull representation of bilinear terms to render a linear OPF model. Suitable pre-processing techniques can be utilized to determine adequate bounds on the nodal voltages and currents to ensure a tight convex relaxation.

\subsection{Problem Formulation}
We denote a general distribution network as a graph $\Gamma(\mathcal{N},\mathcal{E})$, where $\mathcal{N}\coloneqq\{1,...,N\}$ denotes the set of nodes, $\mathcal{E}\coloneqq\{(m,n)\}$ denotes the set of edges, and each phase is expressed as $\phi\in \mathcal{P}, \; \mathcal{P}=\{a,b,c\}$. The CI-OPF for the network is then written as\footnote{For brevity, we omit the McCormick Envelope relaxation of the bilinear terms and the pre-processing step. See \cite{FerroCIM_IFAC,FerroThesis} for details}:
\begingroup\makeatletter\def\f@size{9.5}\check@mathfonts
\def\maketag@@@#1{\hbox{\m@th\normalfont#1}}
\begingroup
\allowdisplaybreaks
\begin{subequations} \label{eq:CI}
\begin{align}
    \underset{x}{\mathop{\min }}\, f & \left(x\right) \label{eq:CIOPF_obj} \\
    AV &= ZI_{\textrm{flow}} \label{eq:CIOPF_ohm} \\
    {{I}^{R}} &= \operatorname{Re}\left( A^TI_{\textrm{flow}} \right),\quad {{I}^{I}} = \operatorname{Im}\left( A^TI_{\textrm{flow}} \right)  \label{eq:CIOPF_Idef} \\
    P_{j}^{\phi} &= V_{j}^{\phi,R}I_{j}^{\phi,R}+V_{j}^{\phi,I}I_{j}^{\phi,I}  \label{eq:CIOPF_Pdef} \\
    Q_{j}^{\phi} &= -V_{j}^{\phi,R}I_{j}^{\phi,I}+V_{j}^{\phi,I}I_{j}^{\phi,R}  \label{eq:CIOPF_Qdef} \\
    \underline{P}_{j}^{\phi} &\le P_{j}^{\phi} \le \overline{P}_{j}^{\phi}, \quad \underline{Q}_{j}^{\phi} \le Q_{j}^{\phi}\le \overline{Q}_{j}^{\phi} \label{eq:CIOPF_powerlim} \\
    \underline{V}_{j}^{\phi,R} &\le V_{j}^{\phi,R} \le \overline{V}_{j}^{\phi,R}, \quad \underline{V}{j}^{\phi,I} \le V_{j}^{\phi,I} \le \overline{V}_{j}^{\phi,I}  \label{eq:CIOPF_Vlim}\\
    \underline{I}_j^{\phi,R} &\le I_j^{\phi,R}\le \overline{I}_j^{\phi,R}, \quad \underline{I}_j^{\phi,I} \le I_j^{\phi,I}\le \overline{I}_j^{\phi,I}  \label{eq:Ilim}
\end{align}
\end{subequations}
\endgroup \endgroup
where $x = \left[{{I}^{R}}\;{{I}^{I}}\;{{V}^{R}}\;{{V}^{I}}\;P\;Q\;I_{\textrm{flow}}^R\;I_{\textrm{flow}}^I\right]$ is the decision vector for the CI-OPF problem; $I, V, P, Q$ denote the vector of nodal current injections, voltages, and real/reactive power injections respectively; $I_{\textrm{flow}}$ denotes the vector of line currents; $A \in \mathbb{R}^{3N\times3N}$ is the 3-phase graph incidence matrix; $Z$ is the system impedence matrix. We use $x^R$ and $x^I$ to denote the real and imaginary components of a complex number $x$; overbar $\overline{{x}}$ and underbar $\underline{{x}}$ denote the upper and lower limits of a variable ${x}$; $\operatorname{Re}(\cdot)$ and $\operatorname{Im}(\cdot)$ denote the real and imaginary components of a complex number. Constraint \eqref{eq:CIOPF_ohm} describes the generalized Ohm's law, \eqref{eq:CIOPF_Idef} describe Kirchhoff's Current Law, and \eqref{eq:CIOPF_Pdef}-\eqref{eq:CIOPF_Qdef} are the definitions of real and reactive power. Constraints \eqref{eq:CIOPF_Pdef}-\eqref{eq:Ilim} are for all nodes $j \in \mathcal{N}$, per each phase $\phi \in \mathcal{P}$, and for all time $t$.
The convex relaxation of the CI-OPF problem through McCormick envelopes adds additional constraints to the problem namely \eqref{eq:Ilim}, concerning the minimum and maximum limits of injected current at each phase of each node. The determination of these upper and lower bounds on nodal current is not trivial. An effective heuristic to define these values is presented in \cite{FerroThesis}.
\subsection{Model of DERs}
We use active sign convention, such that nodal injections are positive. The PV units are modelled as a generator with variable power factor (\textit{pf}), equipped with a multiphase inverter, where the ratio of P and Q determine the \textit{pf} setting:
\begingroup\makeatletter\def\f@size{9.5}\check@mathfonts
\def\maketag@@@#1{\hbox{\m@th\normalfont#1}}
\begin{align}
    P_{j}^{\phi} \tan(\cos^{-1}(-{pf})) &\leq Q_j^{\phi} \leq P_j^{\phi} \tan(\cos^{-1}({pf})) & \label{eq:CIOPF_pf}
\end{align}
\endgroup 

Flexible loads are reductions in real power demand (reactive power loads are fixed), and are modelled as a percentage reduction from the forecasted load. The power constraints are modified as $\overline{P}_j(t) = \underline{P}_j(t)*(1-\alpha_j^\text{DR}(t))$, by recalling the use of active sign convention. For inflexible loads, $\underline{P}_j=\overline{P}_j, \overline{P}_j \leq 0$. Prosumers, which are nodes where both load and generation are present, are modelled with additional variables representing both load and generation, $P^L$ and $P^G$ respectively. For prosumers, the inverter is modelled on $P^G$ and $Q^G$, and load flexibility is modelled on $P^L$. Variables $P^L$ and $P^G$ are both nonnegative, and the same for reactive power.
\vspace{-10pt}
\begingroup\makeatletter\def\f@size{9.5}\check@mathfonts
\def\maketag@@@#1{\hbox{\m@th\normalfont#1}}
\begin{subequations} \label{eq:prosumer}
\begin{align} 
    P_j &= P_j^G-P_j^L, \quad P_j^G \geq 0, P_j^L \geq 0 \\
    Q_j &= Q_j^G-Q_j^L, \quad Q_j^L \geq 0 \\
    \overline{P}_j &= \overline{P}_j^G-\underline{P}_j^L, \quad \underline{P}_j = \underline{P}_j^G-\overline{P}_j^L \\
    \overline{Q}_j &= \overline{Q}_j^G-\underline{Q}_j^L, \quad \underline{Q}_j = \underline{Q}_j^G-\overline{Q}_j^L
\end{align}
\end{subequations}
\endgroup
\noindent Battery storage devices are modelled using the power charge and discharge, $P_{j}^\text{sc}(t)$ and $P_{j}^\text{sd}(t)$ respectively for node $j$ and time $t$. These are nonnegative variables. The state of charge, $b_j(t)$, is calculated as an integral constraint using the actions of the previous period and the initial state of charge, $b_j^0 \coloneqq b_j(t=0)$. We model charge and discharge efficiencies ($\eta_j^C$ and $\eta_j^D$), self-discharge rate ($\eta_j^\text{self}$), and impose a minimum state of charge ($\underline{b}_j$) to ensure battery health. We assume all batteries operate at unity power factor, i.e. $Q_{j} = 0\quad \forall t$.
\begingroup\makeatletter\def\f@size{9.5}\check@mathfonts
\def\maketag@@@#1{\hbox{\m@th\normalfont#1}}
\begin{subequations}
\begin{align}
    P_j &= \frac{1}{\eta_j^\text{D}} P_{j}^\text{sd}(t) - \eta_j^\text{C} P_{j}^\text{sc}(t) \\
    0 &\leq P_{j}^\text{sd} \leq \overline{P}_{j}^\text{sd}, \quad 0 \leq P_{j}^\text{sc} \leq \overline{P}_{j}^\text{sc} \\
    b_{j}(t) &= (1-\eta_j^\text{self})b_{j}(t-1) + \eta_j^\text{C} P_{j}^\text{sc}(t) - \frac{1}{\eta_j^\text{D}} P_{j}^\text{sd}(t) \\
    \underline{b}_{j} &\leq b_{j}(t) \leq \overline{b}_{j}
\end{align}
\end{subequations}
\endgroup
The storage model introduces a dependency on control action in one period to previous periods. The other DERs and loads do not have inter-temporal constraints, so constraints \eqref{eq:CI}-\eqref{eq:prosumer} are simply replicated for each time step $t$.

\section{PAC-Based Distributed Implementation} \label{sec:PAC}
In this section we introduce a PAC-based distributed optimization algorithm to coordinate a large number of spatially distributed DERs at the grids edge. Consider a global optimization problem composed of equality and inequality constraints, which may be coupled in time
\begin{gather}
    \underset{x}{\text{min}} \sum_{i=1}^S f_i(x) \nonumber\\
    \text{s. t.}\quad  Gx=b, \quad Hx\le d\vspace{-0.02in}
    \label{eq:global}
\end{gather}
where $\sum_{{i} \in S} f_i\left({x}\right)$ represents the total objective function. Problem \eqref{eq:global} can be distributed into $K=\left\{1,...,k,...K\right\}$ separate coupled optimization problems, denoted as \textit{atoms}. We use a decomposition profile which separates the vector of variables $x$ into two sets: $\mathsf{L} = \left\{ L_{j}, \quad \forall j \in K\right\}$ and $\mathsf{O} = \left\{ O_{{j}},  \quad \forall j \in K\right\}$, which represent the partition of decision variables ``owned'' and ``copied'' by atom $j$. The set of total variables (owned and copied) by an atom is denoted as $\mathsf{T}$. The decomposition profile also separates the constraints into sets owned by each atom, as $\mathsf{C} = \left\{ C_{{j}}, \quad \forall j \in K \right\}$. The notion of variables copies are used to satisfy the coupling in constraints and/or objective function. In the context of the CI model, the power physics of the grid \eqref{eq:CIOPF_ohm}-\eqref{eq:CIOPF_Idef} result in these coupling constraints. Note that the CI does not have coupling introduced by inequality constraints, however the decomposition can be trivially extended to coupled inequality constraints. Using the decomposition profile, we obtain:
\begin{equation} \label{eq:sof:aso}
	\begin{array}{c}
		\underset {{a}_{{j}}}{\text{{ min}}} \quad \sum_{{j} \in K}{f}_{{j}}\left({a}_{{j}}\right)  \\
		\begin{array}{rl}
			\text{subj. to:} & \left\{\begin{array}{rl}
				{{G}}_{{j}} {a}_{{j}} = {b}_{j}, & \text{ for all } {j} \in K\\
				{{H}}_{{j}} {a}_{{j}} \le {d}_{j}, & \text{ for all } {j} \in K\\
			        B_j a = {0}, & \text{ for all } {j} \in K
			\end{array}\right.
		\end{array}
	\end{array}
\end{equation}
where $a_j$ is atom $j$'s variables (both owned and copied), $f_j \left(a_j\right) $ is the atomic objective function, and $G_j$, $b_j$, $H_j$, and $d_j$ represent the submatrix or subvector of $G, b, H$, and $d$ respectively. Finally, $B$ is in incidence matrix over the owned and copied atomic variables, defined as 
\begin{equation*}
	B_{i}^{m} \triangleq \left\{\begin{array}{rl} -1, & \text{if $i$ is     `owned` and $m$ a related `copy`} \\ 
	    1, & \text{if $m$ is `owned` and $i$ a related `copy`} \\
	    0, & \text{otherwise}\end{array}\right.
\end{equation*}
We then use $B_j$ ($B^j$) to denote the relevant incoming (out-going) edges of the directed graph for atom-$j$. To fully parallelize the optimization, we introduce \textit{coordination constraints}, which must be satisfied for every atom. These require all atomic copied variables in a given ${j}$th atom  to equal the value of their corresponding owned in $i$th atom, $i\neq j$:
\begin{align}
    B_j {a} = 0 \quad \forall {j} \in K
\end{align}

\subsection{An Accelerated Algorithm: NST-PAC}
We next present an accelerated variant of the PAC algorithm in \cite{RomvaryTAC}, which includes time-varying gains and Nesterov type acceleration \cite{nesterov1983method}, called NST-PAC. The NST-PAC is a primal-dual method with $\ell_2$ and proximal regularization, Nesterov type acceleration for both primal and dual variables, and privacy-preserving features. We begin by forming the atomic Lagrangian function:
\begingroup\makeatletter\def\f@size{9.5}\check@mathfonts
\def\maketag@@@#1{\hbox{\m@th\normalfont#1}}
\begin{align}
	\mathcal{L} \left({a},{\mu},{\nu}\right) &= \sum_{{j} \in K} \left[{f}_{{j}} \left({a}_{{j}}\right) + {\mu}_{{j}}^{T} ({{G}}_{{j}} {a}_{{j}}-b_j) + {\nu}_{{j}}^{T}  B_j {a} \right] \nonumber \nonumber \\
		&= \sum_{{j} \in K} \left[{f}_{{j}} \left({a}_{{j}}\right) + {\mu}_{{j}}^{T} ({{G}}_{{j}} {a}_{{j}}-b_j) + {\nu}^{T}  B^j {a}_{j} \right] \nonumber \\
		&\triangleq \sum_{{j} \in K} \mathcal{L}_{{j}} \left({a}_{{j}}, {\mu}_{{j}}, {\nu}\right). \label{eq:pac:alg:1:1}
\end{align}
\endgroup

The algorithm is carried out as below:
\begingroup\makeatletter\def\f@size{9.5}\check@mathfonts
\def\maketag@@@#1{\hbox{\m@th\normalfont#1}}
\begin{align}
    {a}_{{j}} \left[ \tau + 1 \right] &= \underset{{{a}_{{j}} \in \mathbb{R}^{|T_j|}}}{\text{{ argmin}}} \left\{\mathcal{L}_{{j}}\left({a}_{{j}}, \hat{{\mu}}_{{j}} \left[\tau\right], \hat{{\nu}}\left[ \tau \right]\right)\right. 
    +\frac{\rho_j\gamma_j}{2}\left\| {{G}_{j}}{{a}_{j}}-{{b}_{j}} \right\|_{2}^{2}\nonumber\\
    &\quad \left.+ \frac{\rho_j\gamma_j}{2}\left\| {{B}_{j}}{{a}_{j}} \right\|_{2}^{2}+ \frac{1}{2\rho_j} \left\Vert {a}_{{j}} - {a}_{{j}} \left[ \tau \right] \right\Vert_{2}^{2}\right\},\label{pacnst_beg} \\
    {{\hat{a}}_{j}}\left[ \tau +1 \right] &= {{a}_{j}}\left[ \tau +1 \right]+{{\alpha }_{j}}[\tau +1]({{a}_{j}}\left[ \tau +1 \right]-{{a}_{j}}\left[ \tau  \right])\label{pacnst_acc0} \\
    {{\mu }_{j}}\left[ \tau +1 \right] &= {{\hat{\mu }}_{j}}\left[ \tau  \right]+\rho_j\gamma_j({{G}_{j}}{{\hat{a}}_{j}}\left[ \tau +1 \right]-{{b}_{j}}) \\
    {{\hat{\mu }}_{j}}\left[ \tau +1 \right] &= {{\mu }_{j}}\left[ \tau +1 \right]+{{\phi }_{j}}[\tau +1]({{\mu }_{j}}\left[ \tau +1 \right]-{{\mu }_{j}}\left[ \tau  \right])\label{pacnst_acc1} \\
    \text{Comm}&\text{unicate } {{\hat{a}}_{j}}\text{for all }j\in \left[ K \right]\text{ with neighbors} \\
    {{\nu }_{j}}\left[ \tau +1 \right] &= {{\hat{\nu }}_{j}}\left[ \tau  \right] + \rho_j\gamma_j{{B}_{j}}{{\hat{a}}_{j}}\left[ \tau +1 \right] \\
    {{\hat{\nu }}_{j}}\left[ \tau +1 \right] &= {{\nu }_{j}}\left[ \tau +1 \right]+{{\theta }_{j}}[\tau +1]({{\nu }_{j}}\left[ \tau +1 \right]-{{\nu }_{j}}\left[ \tau  \right])\label{pacnst_acc2} \\
    \text{Comm}&\text{unicate } {{\hat{\nu }}_{j}}\text{ for all }j\in \left[ K \right]\text{with}\,\text{neighbors} \label{pacnst_end}
\end{align}
\endgroup

\noindent where $\rho_j$, $\gamma_j$ are atom-varying over-relaxation and step-size parameters, respectively. The proposed NST-PAC uses $\ell2$ regularization terms rather than the prox-linear variant in PAC. Further, both primal and dual variables are accelerated using Nesterov type acceleration, to speed up convergence. Further, we extend the privacy-preserving feature of the PAC algorithm to both the primal and dual variables, by using three iteration-varying and atom specific parameters for the accelerated terms, ${{\alpha }_{j}}\left[ \tau  \right]>\alpha _{j}^{\min }$, ${{\phi }_{j}}\left[ \tau  \right]>\phi _{j}^{\min }$ and ${{\theta }_{j}}\left[ \tau  \right]>\theta _{j}^{\min }$. In the original algorithm privacy is kept only for dual variables. A detailed analysis of convergence rate, communication and computational complexity, and privacy are provided in \cite{RomvaryTAC}.

\section{Case Study} \label{sec:case}
\begin{figure}
    \centering
    \includegraphics[scale=0.4]{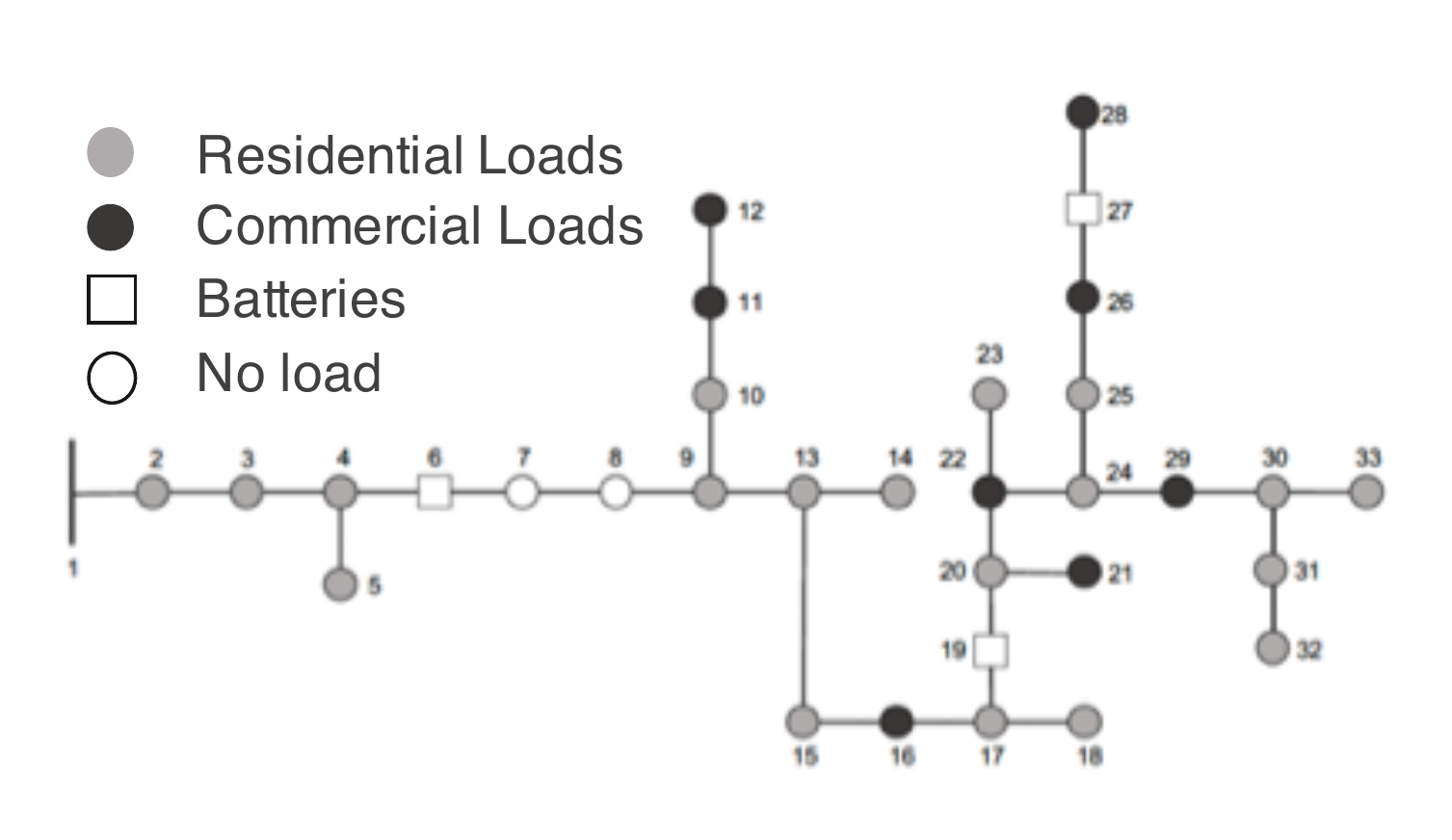}
    \caption{Topology of IEEE-34 network. PV units and flexible loads are present throughout the network.} \vspace{-0.2in}
    \label{fig:ieee_case}
\end{figure}

We consider a case study of San Francisco, California, using the IEEE-34 node network as a proxy for the distribution grid. The load data from the IEEE datasheet serves as the daily average load, and the 24-hr load profiles are obtained from the ODEI dataset from NREL for the Typical Meteorological Year \cite{odei_dataset}. All loads are assumed to have a constant power factor of 0.95. The network loads are classified as residential or commercial loads based on the size of the load, by matching the load levels of the IEEE-34 network with the TMY data. Commercial loads include retail space, small and medium office buildings, primary school, medium and large restaurants, and a hospital. The network is modified to include DERs which include clusters of rooftop photovoltaic (PV) units, flexible loads, and three battery storage units. The penetration of PV is 38\%, as measured by the ratio of nameplate capacity to average system load. This high DER penetration scenario is a reasonable projection given the RPS initiatives in California. Each PV unit is assumed to be equipped with an inverter with corresponding power electronic control, which can be operated at variable power factor in the range of 0.8 to 1. PV curtailment is not considered. The flexible loads are modelled as typical residential cooling loads for California \cite{LBNL_CAISO_DR}. Variations in nodal demand response are obtained by shifting the baseline profile obtained from \cite{LBNL_CAISO_DR} in time and space, with both following zero-mean Gaussian distributions with variances of 0.075 and 0.1 respectively. The three battery units are a 450kW-120kWh community unit\footnote{Modeled on the Ellenbrook unit from the PowerBank trail, Australia \cite{tagabe_2020}} at node 6, a cluster of 40 Tesla Powerwall+ batteries (each 13.5kW-5kWh) at node 19, and a 800kW-185kWh hospital unit at node 27. Any load not met by local generation (or storage) is assumed to be served by the bulk grid at the point of common coupling (PCC), at node 1 in the network ($j = 1$). Figure~\ref{fig:ieee_case} shows the network topology.


\subsection{Scenarios}
We consider three scenarios, as below:
\begin{itemize}
    \item \underline{Scenario A:} \textit{Baseline}. All PV inverters operate at unity pf, and batteries and flexible loads are not present.
    \item \underline{Scenario B:} \textit{Local control}. Each DER owner operates its devices and manages its loads. 
    \item \underline{Scenario C:} \textit{Distributed control}. All devices and loads are coordinated using the NST-PAC algorithm.
\end{itemize}

Scenario A quite trivially is the characteristic duck curve, where the high PV generation at unity power factor results in large ramping requirements of transmission-level generators. Scenario B is a local approach, where each agent will minimize its peak load throughout the day. This serves as a very rough approximation of reducing the ramping requirements of bulk resources, by noticing that the largest ramp typically coincides with the peak demand in the evening. This action can be motivated by the fact that consumers can be charged based on their peak energy consumption. To leverage the capabilities of the storage devices, neighbouring nodes are clustered with the battery. Residential loads (and corresponding DERs) at nodes 3, 4, and 5, share the community battery at node 6. Residential loads at node 20 and the primary school at node 21 share the cluster of Tesla Powerwall+ batteries at node 19. The hospital at node 26 is assumed to own and operate the battery at node 27. All remaining nodes are treated as independent agents. After clustering, there are a total of 26 agents in the network, each managing its own consumption and generation, to minimize its peak load. The multi-period optimization problem solved by each agent is a simple power balance, where any load in excess of local generation is assumed to be served by the bulk system. The power physics between nodes within a cluster is not modelled. The objective function is $f_\text{local}(y) = \text{max}_{t} \left\{ -P(t)\right\}$,
where $y=\left[P, Q,P^L,P^G,Q^L,Q^G,P^\text{sd},P^\text{sc} \right]$.

Scenario C is the PAC-based approach which requires coordination between the devices. This approach accommodates system-level constraints including grid power physics, and minimizes the ramping requirements at the PCC. In Scenario C, the PAC agents solve the CI-based multi-phase unbalanced OPF, to minimize the objective function $f(x)$ which minimizes the difference in power supplied by the PCC from one hour to the next. The function is $f_\text{PAC}(x) = \left\| \sum_{\phi\in \mathcal{P}} P_1^\phi(t) - P_1^\phi(t-1) \right\|$.

\subsection{Results and Discussion}
\begin{figure}[]
    \centering
    \includegraphics[scale=0.28,trim=1.5cm 7cm 0cm 6.5cm,clip]{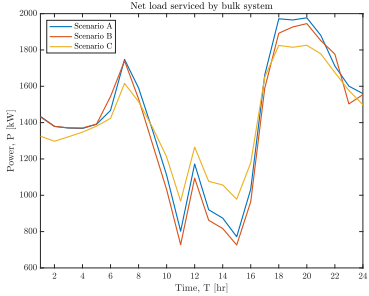}
    \caption{Net load serviced by the bulk system, for each scenario. These plots do not include power loss over lines.}
    \label{fig:res_netload} \vspace{-0.2in}
\end{figure}

\begin{figure}[]
    \centering
    \includegraphics[scale=0.28,trim=1.5cm 7cm 0cm 6.5cm,clip]{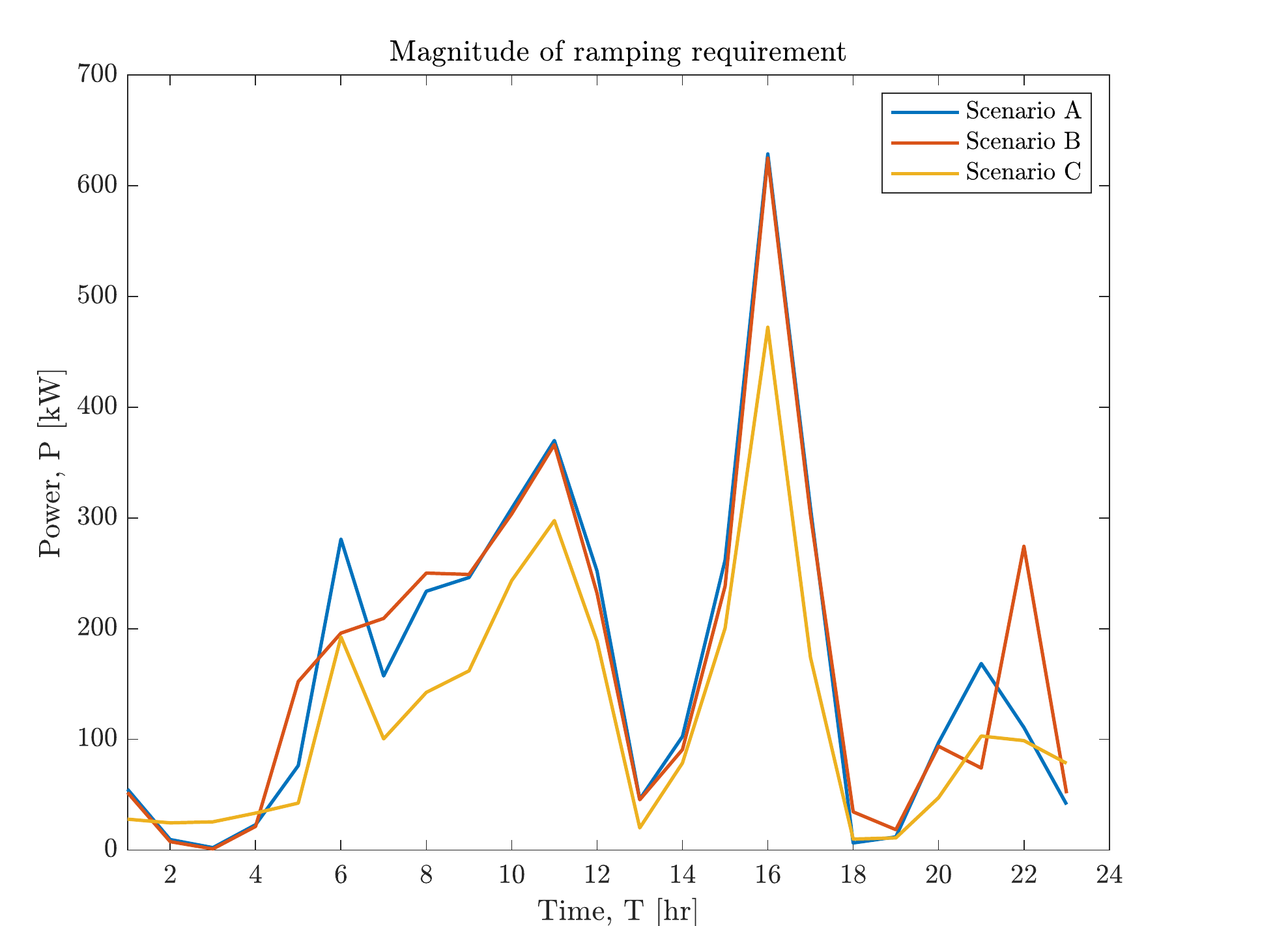}
    \caption{Magnitude of ramping required by bulk system generators for each scenario.}
    \label{fig:res_ramping} \vspace{-0.35in}
\end{figure}

Each scenario was simulated on the IEEE-34 node network. Figure~\ref{fig:res_netload} shows the net load served by the bulk system. The load curve for Scenario A shows the characteristic high ramps down and up when solar generation begins and ends. Scenario B provides minimal improvement, reducing the load for most hours of the day, while Scenario C effectively leverages the DERs to reduce the ramping requirements. Note that the objective function reduces the hour-to-hour change in load, and to do so, increases load during the hours of peak PV generation (roughly 10am to 4pm) to charge the batteries, which are then discharged in the late evening to reduce the net load. The magnitude of hour-to-hour ramping is shown in Figure~\ref{fig:res_ramping}. Notably, the local optimization of Scenario B, which minimizes the peak load throughout the day as a proxy for minimizing the load ramp in the evening, is not able to suitably reduce the ramping requirement. The distributed optimization approach, on the other hand, is able to leverage system-wide information and coordinate the DERs to provide grid-level support, as needed by the bulk system. This coordinated approach is able to reduce ramping requirements throughout the day, with a 23\% reduction in ramping requirements at the 4pm peak. Table~\ref{tab:results} presents the total ramping reduction for Scenarios B and C, as compared to the baseline in A, and the computational run times. As expected, the proposed distributed coordination significantly outperforms the local approach, with 28\% reduction in ramping required. The local approach is unable to provide any reduction. Unsurprisingly, the local approach takes less time to reach a decision (albeit an inferior one), while the distributed approach takes considerably longer (completing 1000 iterations). However, computational time of 20s is still well within the acceptable time-frame for decision making, and enables DERs to provide bulk-level support.

\begin{figure*}[]
	\centering 
	\subfloat[][Minimum initial SOC] {\includegraphics[scale=0.26,trim=1.5cm 7cm 0cm 6.5cm,clip]{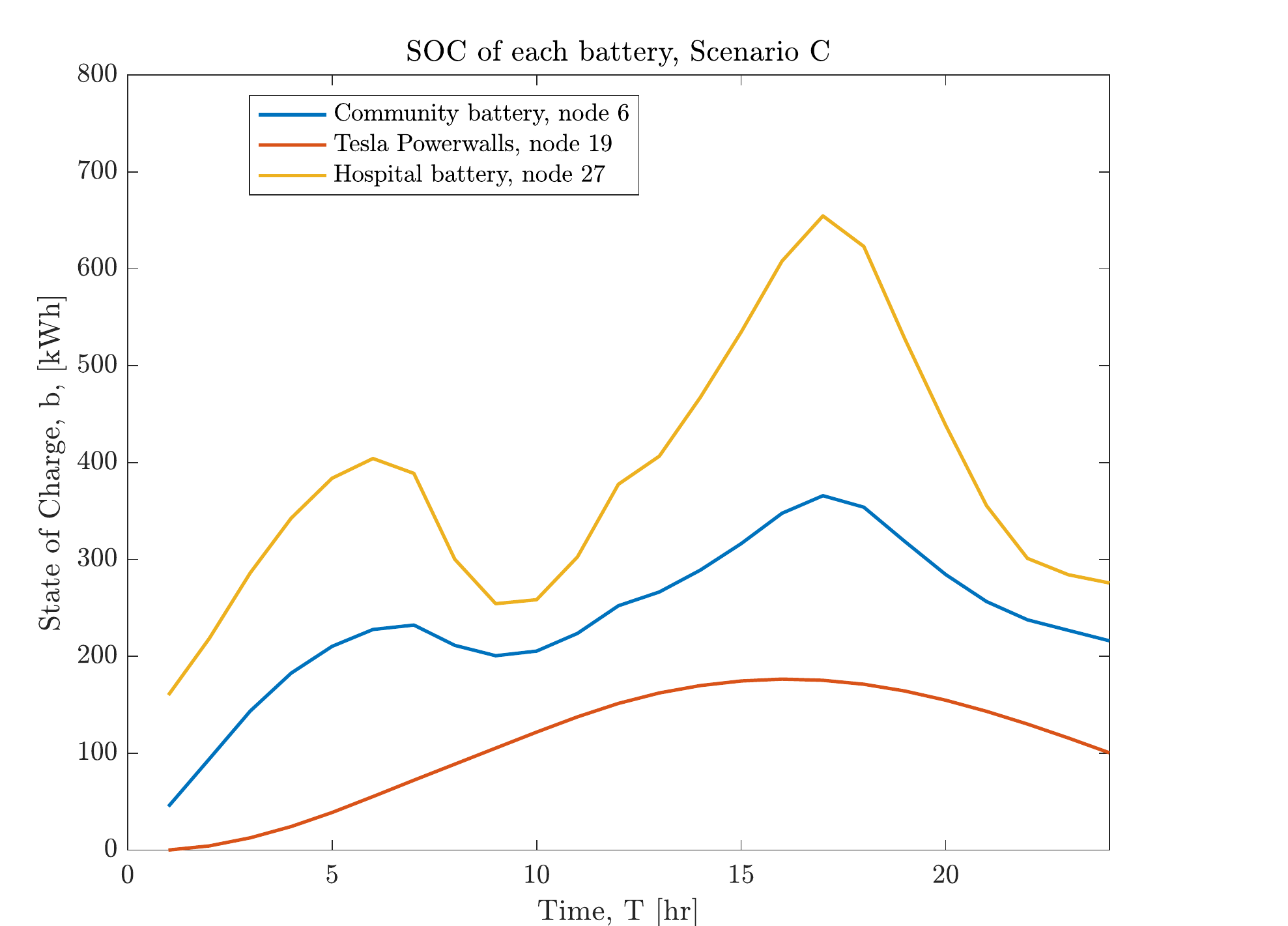}\label{fig:res_storagePAC_zero}}
	\subfloat[][Mid initial SOC] {\includegraphics[scale=0.26,trim=1.5cm 7cm 0cm 6.5cm,clip]{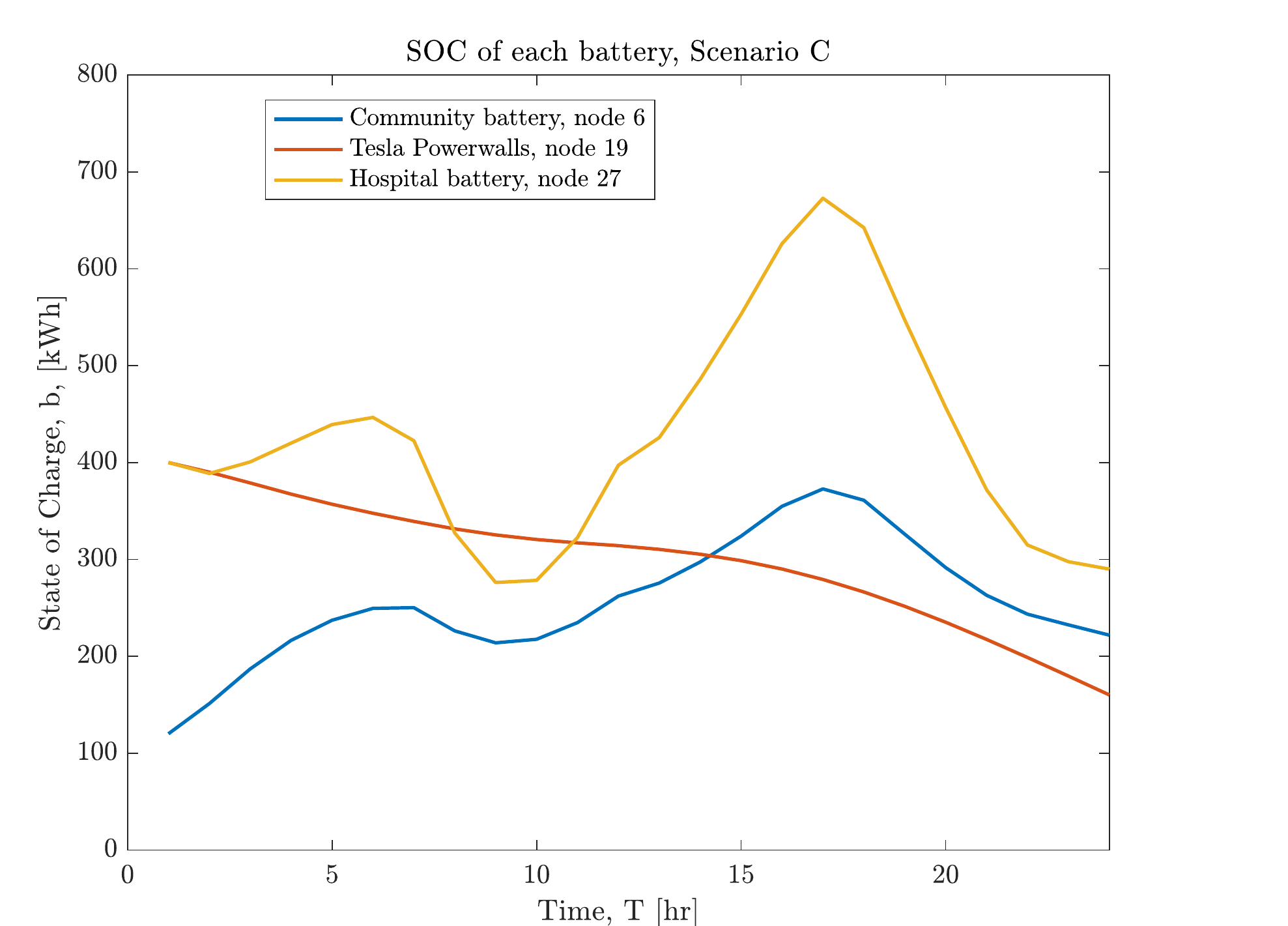}\label{fig:res_storagePAC_mid}}
    \subfloat[][Full initial SOC] {\includegraphics[scale=0.26,trim=1.5cm 7cm 0cm 6.5cm,clip]{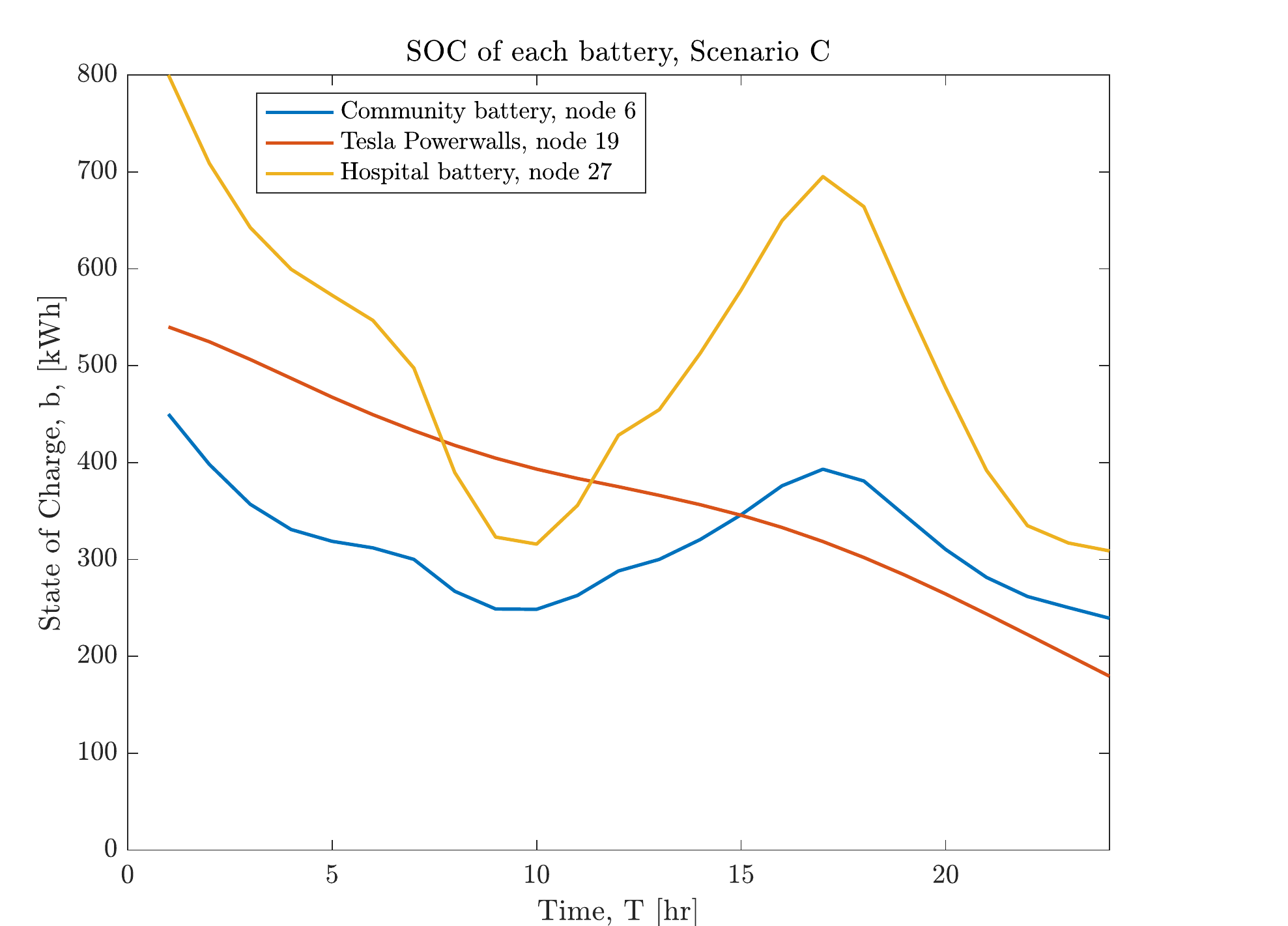}\label{fig:res_storagePAC_full}}
    \caption{State of charge of each storage device for Scenario C, with different initial state of charge for each battery.}
	\label{fig:res_storagePAC} \vspace{-0.2in}
\end{figure*}

We next investigate the impact of the battery units' initial state of charge. We run each Scenario for (1) Minimum SOC where initial capacity is at 45, 0, and 160 kWh; (2) Mid SOC where initial capacity is at 120, 400, and 400 kWh; and (3) Full SOC where initial capacity is at 450, 540, and 800 kWh, respectively for batteries at nodes 6, 19, and 27.
Figure~\ref{fig:res_storagePAC} plots the SOC of each battery unit for the three cases. The usage pattern in each of the cases is very similar, with charging in the early morning and afternoon to build up storage capacity for the evening. Interestingly, the final SOC are non-zero for the three cases, and are quite similar, suggesting the batteries may be able to retain a higher minimum charge to be used as backup power in the case of emergencies. The batteries provide flexibility to increase or decrease load throughout the day, as required by the grid. The plots for net load are quite similar for all three cases, and so have been omitted. 

\begin{table}[]
\begin{tabular}{l|ccc}
 & \multicolumn{1}{l}{A (Baseline)} & \multicolumn{1}{l}{B (Local)} & \multicolumn{1}{l}{C (Proposed)} \\ \hline
Total ramping need (kW) & 3923.7 & 3941.0 & 2839.4 \\
Ramping reduction & - & -0.44\% & 27.63\% \\
Mean run time per agent & - & 0.0947s & 16.97s
\end{tabular}
\caption{Summary of results for Scenarios A thru C. Scenario A is baseline with no decision making, so no ramping reduction or computational time to report.}
\label{tab:results} 
\end{table}

\section{Conclusion} \label{sec:conclusion}
We presented and evaluated a distributed approach to coordinate DERs to provide services to the bulk system. We leverage a CI-based linear model of the unbalanced grid, and an accelerated PAC-based algorithm called NST-PAC. Our case study on a modified IEEE-34 node network shows how distributed techniques can leverage information from different resources to successfully mitigate the duck curve, reducing ramping requirements of bulk system generators by up to 23\%. This framework can be extended to include electric vehicles, by modeling the vehicle-to-grid capabilities and corresponding cost of battery cycling and lifetime degradation in the optimization problem. Using such an approach throughout the distribution grid can reduce the challenges of the new operating condition resulting from high solar and renewable penetration, reduce system costs, and improve renewable integration. Future work will concern the development of faster distributed algorithms applicable to nonconvex costs and constraints, with discrete and continuous controls. 


\bibliographystyle{ieeetr}
\bibliography{references}{}

\end{document}